\documentclass[10pt]{article}

\usepackage{amsmath,amsbsy,amssymb,latexsym,amsthm,dsfont}
\usepackage{a4wide}
\usepackage{graphicx}
\usepackage{cite}
\newtheorem{Theorem}{Theorem}[section]

\newtheorem{Remark}[Theorem]{Remark}
\newtheorem{Proof}{Proof}[section]

\newtheorem{Example}{Example}[section]

\newenvironment{keywords}{\begin{center}
\begin{minipage}[c]{13.4cm} {\bf Keywords:}} {\end{minipage}
\end{center}}

\begin{document}

\title{Numerical solution for fractional variational problems using the Jacobi polynomials}
\author{Hassan Khosravian-Arab$^1$\\
\texttt{h.khosravian@aut.ac.ir}
\and
Ricardo Almeida$^2$\\
\texttt{ricardo.almeida@ua.pt}
}

\date{$^1$ Department of Applied Mathematics, Faculty of Mathematics and Computer Science, Amirkabir University of Technology, No. 424, Hafez Ave. Tehran,Iran\\[0.3cm]
     $^2$ Center for Research and Development in Mathematics and Applications (CIDMA),
Department of Mathematics, University of Aveiro, 3810--193 Aveiro, Portugal}

\maketitle

\begin{abstract}
We exhibit a numerical method to solve fractional variational problems, applying a decomposition formula based on Jacobi polynomials. Formulas for the fractional derivative and fractional integral  of the Jacobi polynomials are proven. By some examples, we show the convergence of such procedure, comparing the exact solution with numerical approximations.
\end{abstract}


\section{Fractional variational calculus}

Variational calculus deals with optimization problems for
functionals depending on some variable function $y$ and some
derivative of $y$  (see e.g. \cite{Gelfand,Sagan}). In many
cases, the dynamic of such trajectories are not described by
integer-order derivatives, but by real-order derivatives
\cite{Kilbas,Podlubny}. Solving these kind of problems usually
implies finding the solutions of a fractional differential
equation, the so-called Euler-Lagrange equation
\cite{Agrawal,Atan1,Baleanu,Jumarie,Malinowska0,Malinowska}. The main problem that arises with this approach is that in most cases there is no way to determine the exact solution. To overcome this situation, many numerical methods
are being developed at this moment for fractional problems. One
of the more commonly used methods consists in approximating the
function by a polynomial $y_n$ and the fractional derivative of
$y$ by the fractional derivative of $y_n$, and by doing this we
rewrite the initial problem in a way such that applying already
known methods  from numerical analysis we can determine the
optimal solution.

The variational problem that we address in this paper is stated in the following way.
Given $\alpha,\beta\in(0,1)$, determine the minimizers of
\begin{equation}\label{clasic}J[y]=\int_a^b F(x,y(x),{_aD_x^{\alpha}}y(x),{_aI_x^{\beta}}y(x))dx,\end{equation}
under the constraint
\begin{equation}\label{bcond0}{_aI_x^{1-\alpha}}y(x)\Big|_{x=b}=y_b.\end{equation}
Here, ${_aD_x^{\alpha}}y(x)$ denotes the Riemann-Liouville fractional derivative of $y$ of order $\alpha$,
$${_aD_x^{\alpha}}y(x)=\frac{1}{\Gamma(1-\alpha)}\frac{d}{dx}\int_a^x (x-t)^{-\alpha}y(t)dt,$$
and ${_aI_x^{\beta}}y(x)$ the Riemann-Liouville fractional integral of $y$ of order $\beta$,
$${_aI_x^{\beta}}y(x)= \frac{1}{\Gamma(\beta)}\int_a^x (x-t)^{\beta-1}y(t) dt.$$
We note that since function $y$ is continuous, the condition
$\displaystyle {_aI_x^{1-\alpha}}y(x)\Big|_{x=a}=0$ appears
implicitly.

\section{Numerical Method}
In this section we present a numerical method to solve the
problem presented in Eqs. \eqref{clasic}-\eqref{bcond0}. We can find several methods in the literature to  solve fractional problem types \cite{Atanackovic1,Pooseh0,Zhuang}. Our main
idea is described in the following way: by using the Jacobi
polynomials, the initial problem is converted into a non-linear
programming problem, without dependence of fractional derivatives
and fractional integrals. By doing this, we are able to find an
approximation for the minimizers of the functional. To start, we
briefly review some basic definitions of Jacobi polynomials.

\subsection{Jacobi polynomials}
The Jacobi polynomials $P_{n}^{(\alpha,\beta)}(t)$ of indices $\alpha,\ \beta$ and degree $n$ are defined by
\begin{equation}\label{jaco}
P_n^{(\alpha,\beta)}(t)=\sum_{k=0}^{n}\frac{(-1)^{n-k}(1+\beta)_{n}(1+\alpha+\beta)_{n+k}}{k!(n-k)!(1+\beta)_{k}(1+\beta+\alpha)_{n}}\left(\frac{t+1}{2}\right)^k,
\end{equation}
where $\alpha,\beta>-1$ are real parameters and
\[
(a)_0=1,\ \ (a)_{i}=a(a+1)\ldots(a+i-1),
\]
The Jacobi polynomials are mutually orthogonal over the interval $(-1, 1)$ with respect to the weight function $w^{\alpha,\beta}(t)=(1-t)^{\alpha}(1+t)^{\beta}$. The Jacobi polynomials $P_{n}^{(\alpha,\beta)}(t)$ reduce to the Legendre polynomials $P_n(t)$ for $\alpha=\beta=0$, and to the Chebyshev polynomials $T_n(t)$ and $U_n(t)$ for $\alpha = \beta = \mp1/2$, respectively \cite{Canu:07:Spe}.

Another useful definition of the Jacobi polynomials of indices $\alpha,\ \beta$ and degree $n$ is as  \cite{Canu:07:Spe,Gaut:04:Ort}:
\begin{equation}\label{jacoo}
P_n^{(\alpha,\beta)}(t)=\frac{(-1)^n}{2^nn!}(1-t)^{-\alpha}(1+t)^{-\beta}\frac{d^n}{dt^n}[(1-t)^{\alpha+n}(1+t)^{\beta+n}].
\end{equation}
This is a direct generalization of the Rodrigues formula
for the Legendre polynomials, to which it reduces for
$\alpha=\beta = 0$.

To present our numerical method, we use three interesting
theorems as follows.

\begin{Theorem}\label{mmm} \cite{esmsha:11:pse} Let $\alpha>0$ be
a real number and $x \in [a,b]$. Then,
$${_aD_x^{\alpha}} [ (x-a)^\alpha P_k^{(0,\alpha)}(\frac{2(x-a)}{b-a}-1)]=\frac{\Gamma(k+\alpha+1)}{\Gamma(k+1)} P_k^{(\alpha,0)}(\frac{2(x-a)}{b-a}-1).$$
\end{Theorem}

\begin{Theorem}\label{our} Let $\alpha-\beta>-1,\ \beta>-1$ be two real numbers and $x \in [a,b]$. Then,
$${_aD_x^{\beta}} [ (x-a)^\alpha P_k^{(0,\alpha)}(\frac{2(x-a)}{b-a}-1)]=\frac{\Gamma(k+\alpha+1)}{\Gamma(k+\alpha-\beta+1)}(x-a)^{\alpha-\beta} P_k^{(\beta,\alpha-\beta)}(\frac{2(x-a)}{b-a}-1).$$
\end{Theorem}
\begin{Proof}
By substituting $t=\frac{2(x-a)}{b-a}-1$ in (\ref{jaco}), we get
\begin{equation}\label{eqt3}
\xi(x):=(x-a)^{\alpha} P_{k}^{(0,\alpha)}(\frac{2(x-a)}{b-a}-1)
=\sum_{m=0}^{k}\frac{(-1)^{k-m}(1+\alpha)_{k+m}}{m!(k-m)!(1+\alpha)_m}\frac{(x-a)^{m+\alpha}}{(b-a)^m}
\end{equation}
Taking the Riemann-Liouville fractional derivative of order $\alpha$ on both side of (\ref{eqt3}), we conclude
\begin{eqnarray*}
{}_aD_x^{\beta}\xi(x)&=&\sum_{m=0}^{k}\frac{(-1)^{k-m}(1+\alpha)_{k+m}\Gamma(m+\alpha+1)}{m!(k-m)!(1+\alpha)_{m}\Gamma(m+\alpha-\beta+1)}\frac{(x-a)^{m+\alpha-\beta}}{(b-a)^m}\nonumber\\
&=&\frac{(1+\alpha)_k\Gamma(\alpha+1)(x-a)^{\alpha-\beta}}{(1+\alpha-\beta)_k\Gamma(\alpha-\beta+1)}\sum_{m=0}^{k}\frac{(-1)^{k-m}(1+\alpha-\beta)_k(1+\alpha)_{k+m}}{m!(k-m)!(1+\alpha-\beta)_m(1+\alpha)_k}(\frac{x-a}{b-a})^m\nonumber\\
&=&\frac{\Gamma(k+\alpha+1)}{\Gamma(k+\alpha-\beta+1)}(x-a)^{\alpha-\beta} P_k^{(\beta,\alpha-\beta)}(\frac{2(x-a)}{b-a}-1),
\end{eqnarray*}
and the proof is completed.
\end{Proof}
We recall that Theorem \ref{our} is a generalized form of the
Theorem \ref{mmm} that was proved in \cite{esmsha:11:pse}.
\begin{Theorem}\label{ouri} Let $\alpha+\beta>-1,\ \beta<1$ be two real numbers and $x \in [a,b]$. Then,
$${_aI_x^{\beta}} [ (x-a)^\alpha P_k^{(0,\alpha)}(\frac{2(x-a)}{b-a}-1)]=\frac{\Gamma(k+\alpha+1)}{\Gamma(k+\alpha+\beta+1)}(x-a)^{\alpha+\beta} P_k^{(-\beta,\alpha+\beta)}(\frac{2(x-a)}{b-a}-1).$$
\end{Theorem}
\begin{Proof}
This theorem is proved if we replace $\beta$ by $-\beta$ in Theorem \ref{our}.
\end{Proof}
\subsection{Presented method}

Our aim is to solve the following variational problem:
\begin{equation}\label{Lei}
J[y]=\int_a^b F(x,y(x),{}_{a}D_{x}^{\alpha}y(x),{}_{a}I_{x}^{\beta}y(x))dx \quad \rightarrow \quad \min,
\end{equation}
under the constraint
\begin{equation}\label{bc}
{}_{a}I_{x}^{1-\alpha}y(x)\Big|_{x=b}=y_b.
\end{equation}
We remark that when $\alpha=1$ and $\beta=0$, we recover the fundamental problem:
$$\int_a^b F(x,y(x),y'(x))dx \quad \rightarrow \quad \min,$$
under the constraints
$$y(a)=0 \quad \mbox{and} \quad y(b)=y_b.$$
To solve the problem, we approximate $y(x)$ by the formula
\begin{eqnarray}\label{type1}
y(x)\approx{y}_{n}(x)&=&\sum_{i=0}^n c_i (x-a)^\alpha \left (P_i^{(0,\alpha)}(\frac{2(x-a)}{b-a}-1)\right),
\end{eqnarray}
where $c_i,\ i=0,1,2,\ldots,n$ are unknown coefficients that should be determined.

Using Theorems \ref{mmm} and \ref{ouri}, we can obtain ${}_{a}D_{x}^{\alpha} {y}_n(x)$ and ${}_{a}I_{x}^{\beta} {y}_n(x)$  as:
\begin{eqnarray}\label{Refo1}
{_aD_x^{\alpha}}y(x)\approx{_aD_x^{\alpha}}{y}_{n}(x)&=&\sum_{i=0}^n c_i\left (\frac{\Gamma(i+\alpha+1)}{\Gamma(i+1)}P_i^{(\alpha,0)}(\frac{2(x-a)}{b-a}-1)\right),
\end{eqnarray}
\begin{eqnarray}\label{Refo2}
{_aI_x^{\beta}}y(x)\approx{_aI_x^{\beta}}{y}_{n}(x)&=&\sum_{i=0}^n c_i(x-a)^{\alpha+\beta}\left (\frac{\Gamma(i+\alpha+1)}{\Gamma(i+\alpha+\beta+1)}P_i^{(-\beta,\alpha+\beta)}(\frac{2(x-a)}{b-a}-1)\right).\nonumber\\
\end{eqnarray}

By substituting \eqref{type1}-\eqref{Refo2} in $J$ and using a quadrature rule, we can approximate $J(y)$ as:

\begin{eqnarray}
  J(y)\approx J_{n}(y)&=&\int_a^b F\left(x,y_n(x),{_aD_x^{\alpha}}y_n(x),{_aI_x^{\beta}}y_n(x)\right)dx\nonumber\\
                    &\approx&\sum_{j=0}^k \omega_j F\left(\xi_j,y_n(\xi_j),{_aD_x^{\alpha}}y_n(\xi_j),{_aI_x^{\beta}}y_n(\xi_j)\right),\label{Z1}
\end{eqnarray}
subject to:
\begin{eqnarray}\label{Refo3}
{_aI_x^{1-\alpha}}y(x)\Big|_{x=b}&\approx&{_aI_x^{1-\alpha}}{y}_{n}(x)\Big|_{x=b}\nonumber\\
&=&\sum_{i=0}^n c_i(b-a)\left (\frac{\Gamma(i+\alpha+1)}{\Gamma(i+2)}P_i^{(\alpha-1,1)}(1)\right)=y_b,
\end{eqnarray}
 where $\xi_j$ and $\omega_j$ are the nodes and weights of quadrature rule. In order to obtain a high order accuracy, we use the Gauss-Legendre quadrature rule\cite{Gaut:04:Ort,Tre:00:Spe}.  Note that the above approximation can be considered as a function of the unknown parameters $c_0,c_1,\cdots,c_n$.
\\

Finally, the problem \eqref{Lei}-\eqref{bc} is converted to a mathematical programming problem with the unknown parameters $c_0,c_1,\cdots,c_n$, as:
$$I(c_0,c_1,\cdots,c_n)=\sum_{j=0}^k \omega_j F\left(\xi_j,y_n(\xi_j),{_aD_x^{\alpha}}y_n(\xi_j),{_aI_x^{\beta}}y_n(\xi_j)\right) \quad \rightarrow \quad \min, $$
subject to
$$\sum_{i=0}^n c_i(b-a)\left (\frac{\Gamma(i+\alpha+1)}{\Gamma(i+2)}P_i^{(\alpha-1,1)}(1)\right)=y_b.$$


\section{Numerical results}

To test the efficiency of the procedure, we will study a fractional variational problem with known solution, and after we compare it with some numerical solutions. The procedure of the following theorem is based on the fractional Leitmann's principle, as showed in \cite{AlTor:10:Lei}.

\begin{Theorem}\label{TeoExample}
Let $g$ and $h$ be two functions of class $\mathcal{C}^1$ with $g(x)\neq0$ on $[a,b]$, and $\beta$ and $\epsilon$ real numbers with $\beta\in(0,1)$. The global minimizer of the fractional variational problem
\begin{equation}\label{Prob1}
J[y]=\int_a^b \left(g(x){}_{a}D_{x}^{1-\beta}y(x)+g'(x){}_aI_x^{\beta}y(x)+h'(x)\right)^2 dx \quad \rightarrow \quad \min,
\end{equation}
under the constraint
\begin{equation}\label{bcond}
{}_{a}I_{x}^{\beta}y(x)\Big|_{x=b}=\epsilon,
\end{equation}
is given by the function
\begin{equation}\label{exact}
y(x)={}_{a}D_{x}^{\beta}\left[ \frac{Ax+C-h(x)}{g(x)}\right],
\end{equation}
where
\[
A=\frac{g(b)\epsilon+h(b)-h(a)}{b-a},\ C=\frac{bh(a)-ah(b)-ag(b)\epsilon}{b-a}.
\]
\end{Theorem}
\begin{Proof}
We know that
\[
\frac{d}{dx}\left[g(x){}_aI_x^{\beta}y(x)+h(x)\right]=g(x){}_{a}D_{x}^{1-\beta}y(x)+g'(x){}_aI_x^{\beta}y(x)+h'(x).
\]
Consider the transformation $y(x)=\tilde{y}(x)+f(x)$, where
$y(x)$ is a function that satisfies problem
\eqref{Prob1}-\eqref{bcond} and $f(x)$ is an unknown function
that to be determined later. Then
\begin{multline*}
\left(\frac{d}{dx}\left[g(x){}_a I_x^{\beta}\tilde{y}(x)+h(x)+g(x){}_a I_x^{\beta}f(x)\right]\right)^2-
\left(\frac{d}{dx}\left[g(x){}_a I_x^{\beta}\tilde{y}(x)+h(x)\right]\right)^2 \\
=2 \frac{d}{dx}\left[g(x){}_a I_x^{\beta}f(x)\right]\frac{d}{dx}\left[g(x){}_a I_x^{\beta}\tilde{y}(x)+h(x)\right]+\left(\frac{d}{dx}\left[g(x){}_a I_x^{\beta}f(x)\right]\right)^2 \\
= \frac{d}{dx}\left[g(x){}_a I_x^{\beta}f(x)\right]\frac{d}{dx}\left[2g(x){}_a I_x^{\beta}\tilde{y}(x)+2h(x)+g(x){}_a I_x^{\beta}f(x)\right].\qquad\qquad
\end{multline*}
Let $f$ be such that
\[
\frac{d}{dx}\left[g(x){}_a I_x^{\beta}f(x)\right]=const.
\]
Integrating, we get
\[
f(x)={}_{a}D_{x}^{\beta}\left(\frac{Ax+B}{g(x)}\right).
\]
Now, consider the new problem
\begin{equation}\label{Prob11}
J[y]=\int_a^b \left(g(x){}_{a}D_{x}^{1-\beta}\tilde{y}(x)+g'(x){}_aI_x^{\beta}\tilde{y}(x)+h'(x)\right)^2 dx \quad \rightarrow \quad \min,
\end{equation}
under the constraint
\begin{equation}\label{bcond1}
{}_{a}I_{x}^{\beta}\tilde{y}(x)\Big|_{x=b}=\frac{1-h(b)}{g(b)}.
\end{equation}
It is easy to see that
\[
\tilde{y}(x)={}_{a}D_{x}^{\beta}\left(\frac{1-h(x)}{g(x)}\right),
\]
is a solution of problem \eqref{Prob11}-\eqref{bcond1}. Therefore,
\[
y(x)={}_{a}D_{x}^{\beta}\left(\frac{Ax+C-h(x)}{g(x)}\right),\
C=B+1,
\]
is a solution of problem \eqref{Prob1}-\eqref{bcond}. Using the boundary conditions
\[
{}_{a}I_{x}^{\beta}y(x)\Big|_{x=a}=0,\ {}_{a}I_{x}^{\beta}y(x)\Big|_{x=b}=\epsilon,
\]
we obtain the values of constants $A$ and $C$ as:
\[
A=\frac{g(b)\epsilon+h(b)-h(a)}{b-a},\ C=\frac{bh(a)-ah(b)-ag(b)\epsilon}{b-a}.
\]
\end{Proof}
\begin{Remark}
For $\beta=0$ the problem \eqref{Prob1}-\eqref{bcond} coincides with the classical problem of the calculus of variations
\begin{equation}
J[y]=\int_a^b \left(g(x)y'(x)+g'(x)y(x)+h'(x)\right)^2 dx \quad \rightarrow \quad \min,
\end{equation}
under the constraint
\begin{equation}
y(a)=0,\ y(b)=\epsilon.
\end{equation}
In this case the global minimizer is obtained from \eqref{exact} as:
\begin{equation*}
y(x)= \frac{Ax+C-h(x)}{g(x)},
\end{equation*}
where
\[
A=\frac{g(b)\epsilon+h(b)-h(a)}{b-a},\ C=\frac{bh(a)-ah(b)-ag(b)\epsilon}{b-a}.
\]
\end{Remark}
\begin{Remark}
For $\beta=1-\alpha$, $\alpha\in(0,1)$ and $g(x)=h(x)$ the problem \eqref{Prob1}-\eqref{bcond} reduces to the problem
\begin{equation}\label{Prob2}
J[y]=\int_a^b \left(g(x){}_{a}D_{x}^{\alpha}{y}(x)+g'(x)({}_aI_x^{1-\alpha}y(x)+1)\right)^2 dx \quad \rightarrow \quad \min,
\end{equation}
under the constraint
\begin{equation}\label{bcond2}
{}_{a}I_{x}^{1-\alpha}y(x)\Big|_{x=b}=\epsilon.
\end{equation}
In this case the global minimizer is obtained from \eqref{exact} as:
\[
y(x)={}_{a}D_{x}^{1-\alpha}\left(\frac{Ax+C}{g(x)}-1\right),
\]
where
\[
A=\frac{g(b)(\epsilon+1)-g(a)}{b-a},\ C=\frac{bg(a)-ag(b)(\epsilon+1)}{b-a}.
\]

This problem for $a=0,\ b=1$ was studied in \cite{AlTor:10:Lei}.
\end{Remark}
\begin{Remark}
For $g(x)=1$, $h(x)=0$ and $\beta=1-\alpha$, the problem \eqref{Prob2}-\eqref{bcond2} reduces to
\begin{equation}\label{Prob3}
J[y]=\int_a^b \left({}_{a}D_{x}^{\alpha}{y}(x)\right)^2 dx \quad \rightarrow \quad \min,
\end{equation}
under the constraint
\begin{equation}\label{bcond3}
{}_{a}I_{x}^{1-\alpha}y(x)\Big|_{x=b}=\epsilon.
\end{equation}
In this case the global minimizer is obtained from \eqref{exact} as:
\[
y(x)=\frac{\epsilon}{b-a}\left(\frac{x^\alpha}{\Gamma(1+\alpha)}-a\frac{x^{\alpha-1}}{\Gamma(\alpha)}\right).
\]

This problem for $a=0,\ b=1$ was studied in \cite{AlTor:10:Lei}.
\end{Remark}
\begin{Remark}
For $\beta=0$ and $g(x)=h(x)$ the problem \eqref{Prob2}-\eqref{bcond2} coincides with the classical problem of the calculus of variations
\begin{equation}
J[y]=\int_a^b \left(g(x)y'(x)+g'(x)(y(x)+1)\right)^2 dx \quad \rightarrow \quad \min,
\end{equation}
under the constraint
\begin{equation}
y(a)=0,\ y(b)=\epsilon.
\end{equation}
In this case the global minimizer is obtained from \eqref{exact} as:
\begin{equation*}
y(x)= \frac{Ax+C}{g(x)}-1,
\end{equation*}
where
\[
A=\frac{g(b)(\epsilon+1)-g(a)}{b-a},\ C=\frac{bg(a)-ag(b)(\epsilon+1)}{b-a}.
\]
We note that this problem for $a=0,\ b=1$ has been studied by Leitmann in \cite{Lei:67:Ano}.
\end{Remark}

\begin{figure}
  \center\includegraphics[width=7cm]{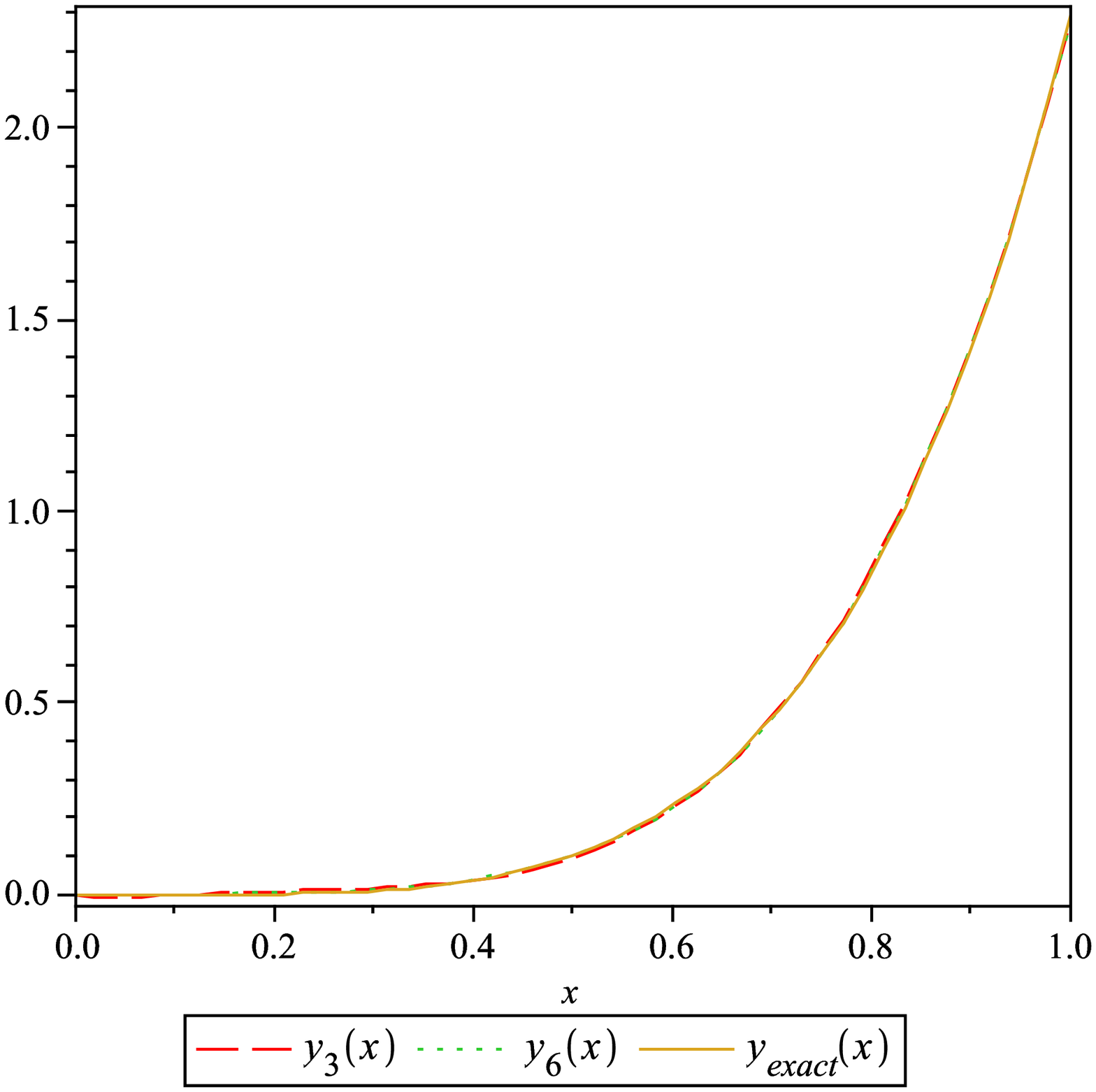}\includegraphics[width=7cm]{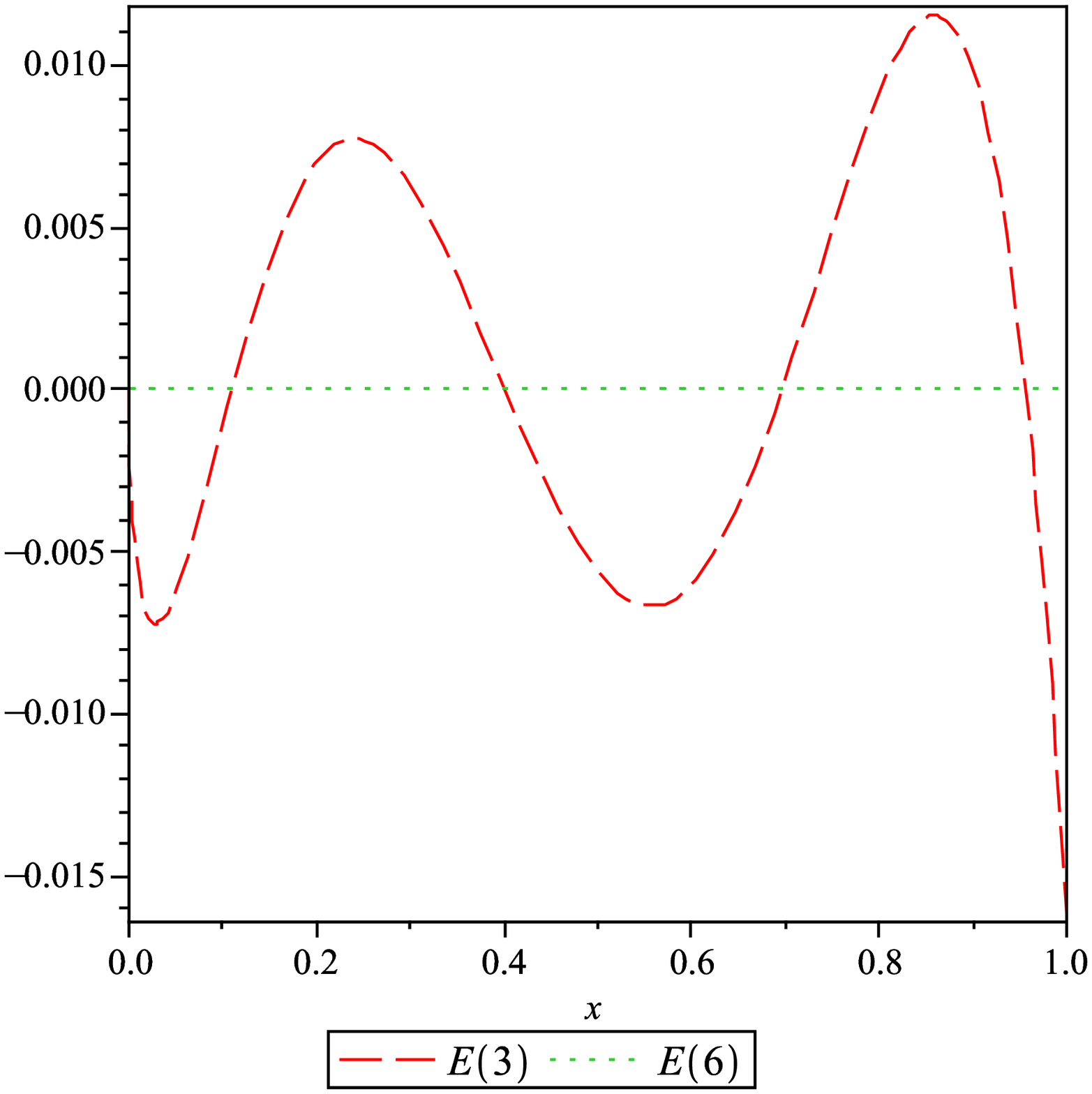}
  \caption{Comparison of exact solution and numerical solution for $n=3,\ 6$ and $\alpha=0.5,\ \beta=5,\ \epsilon=1$ (left) and error between exact solution and numerical solution $E(n)={y}_n(x)-y_{exact}(x)$ for $n=3,\ 6$ and $\alpha=0.5,\ \beta=5,\ \epsilon=1$ (right) in Example \ref{ex1}.}\label{Fig.1.}
\end{figure}

\begin{figure}
  \center\includegraphics[width=7cm]{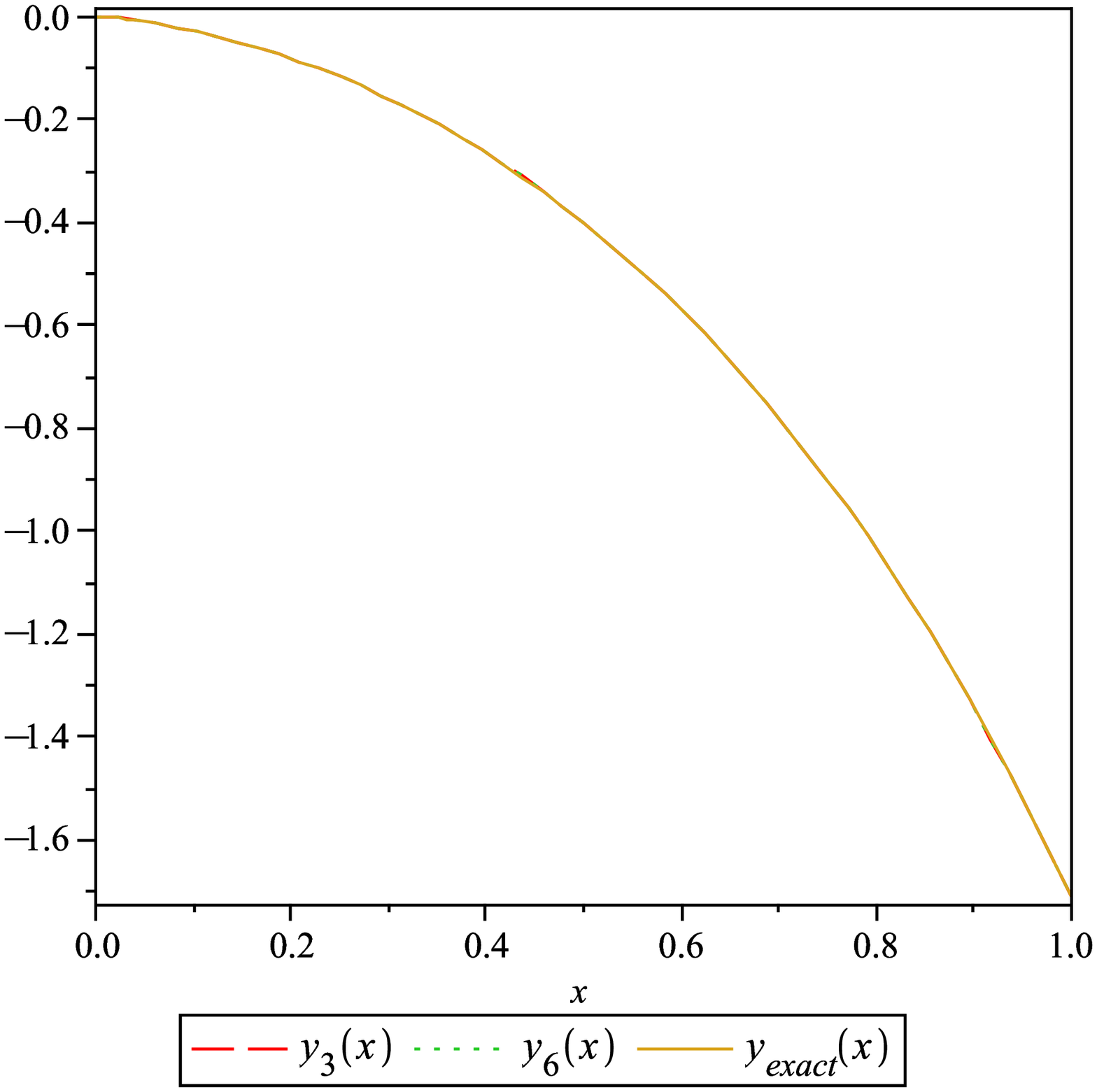}\includegraphics[width=7cm]{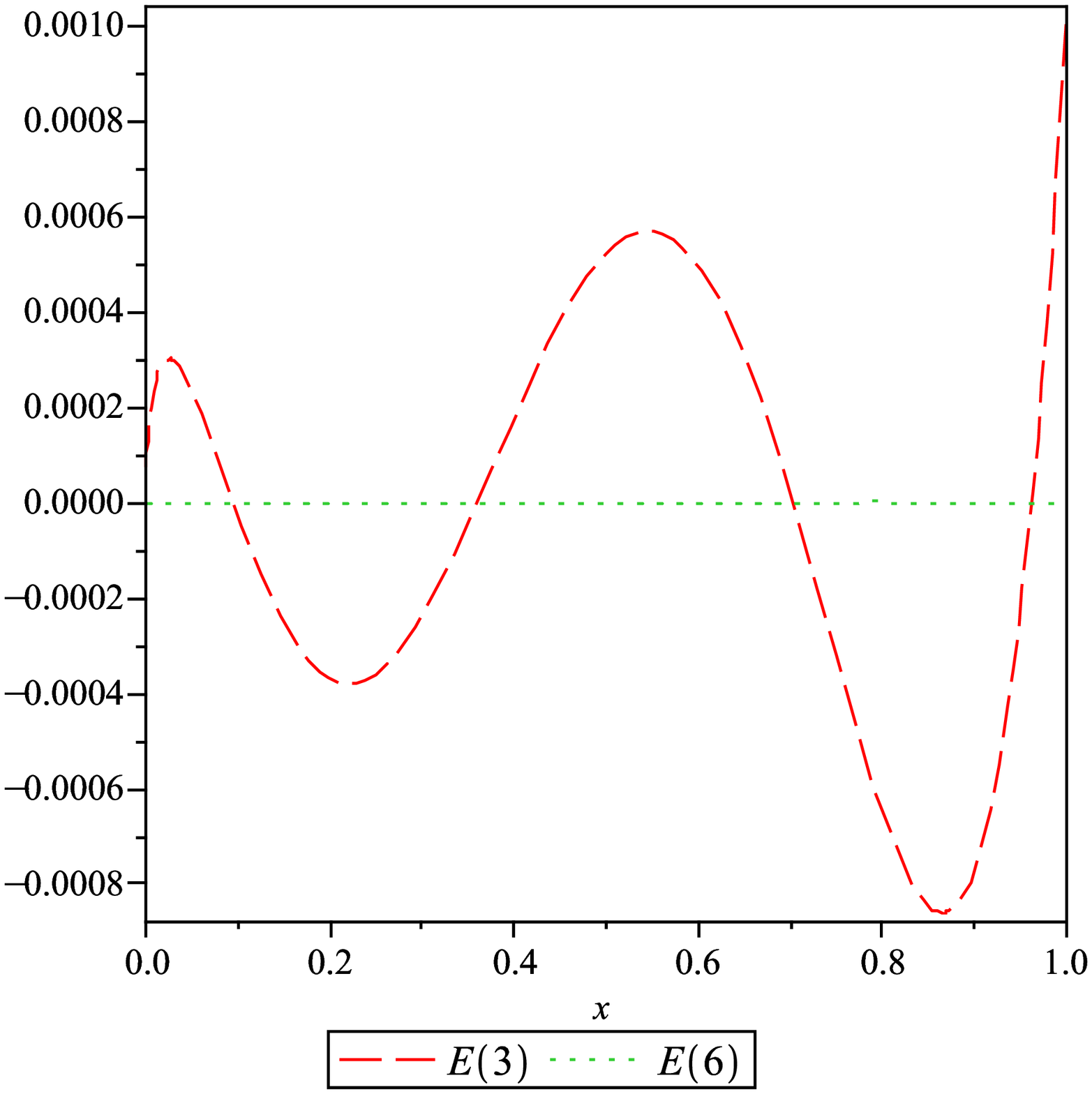}
  \caption{Comparison of exact solution and numerical solution for $n=3,\ 6$ and $\alpha=0.5$ and $\nu=1,\ \epsilon=-1$ (left) and error between exact solution and numerical solution $E(n)={y}_n(x)-y_{exact}(x)$ for $n=3,\ 6$ and $\alpha=0.5$ and $\nu=1,\ \epsilon=-1$ (right) in Example \ref{ex2}.}\label{Fig.2.}
\end{figure}

\begin{figure}
  \center\includegraphics[width=7cm]{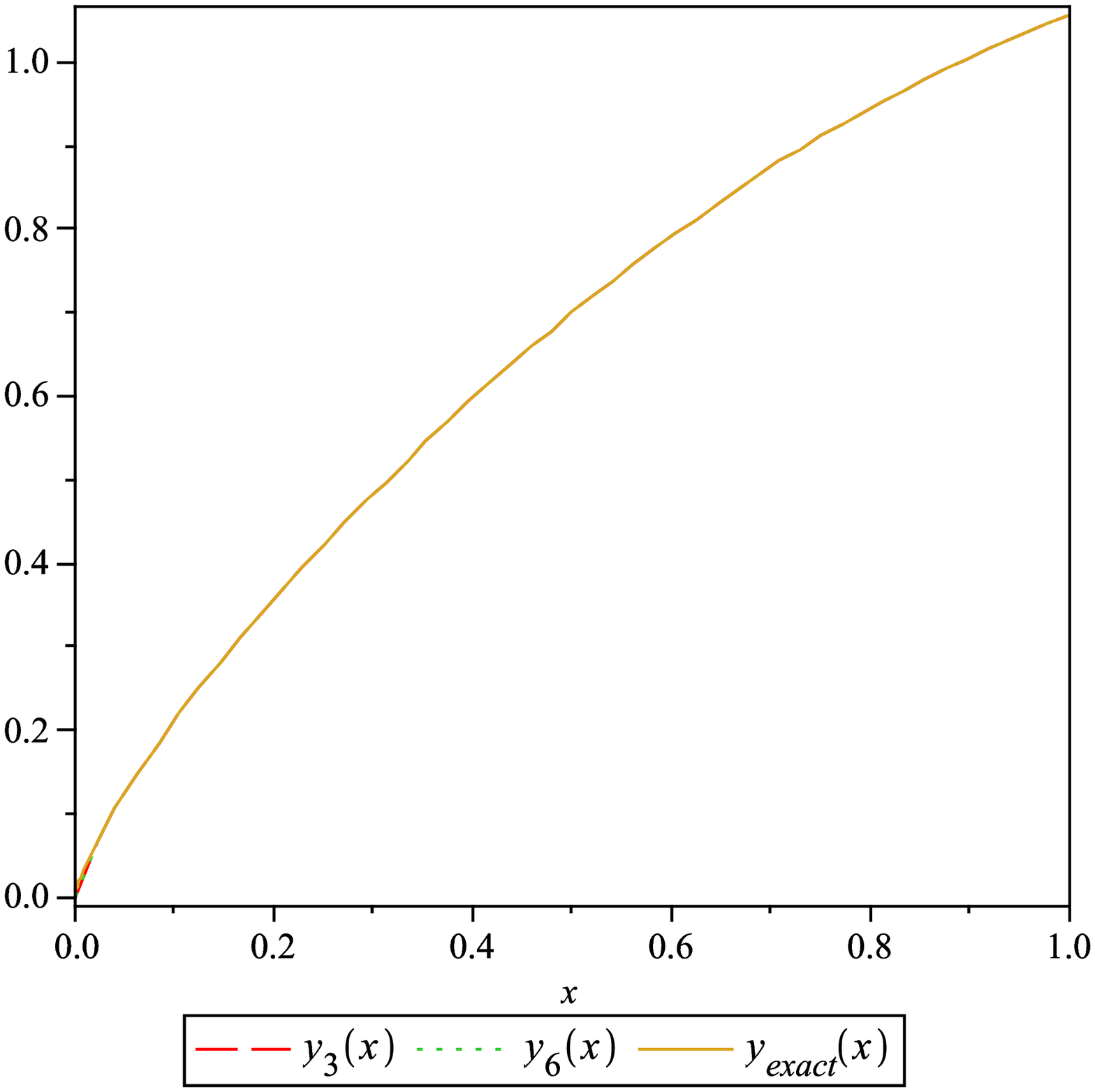}\includegraphics[width=7cm]{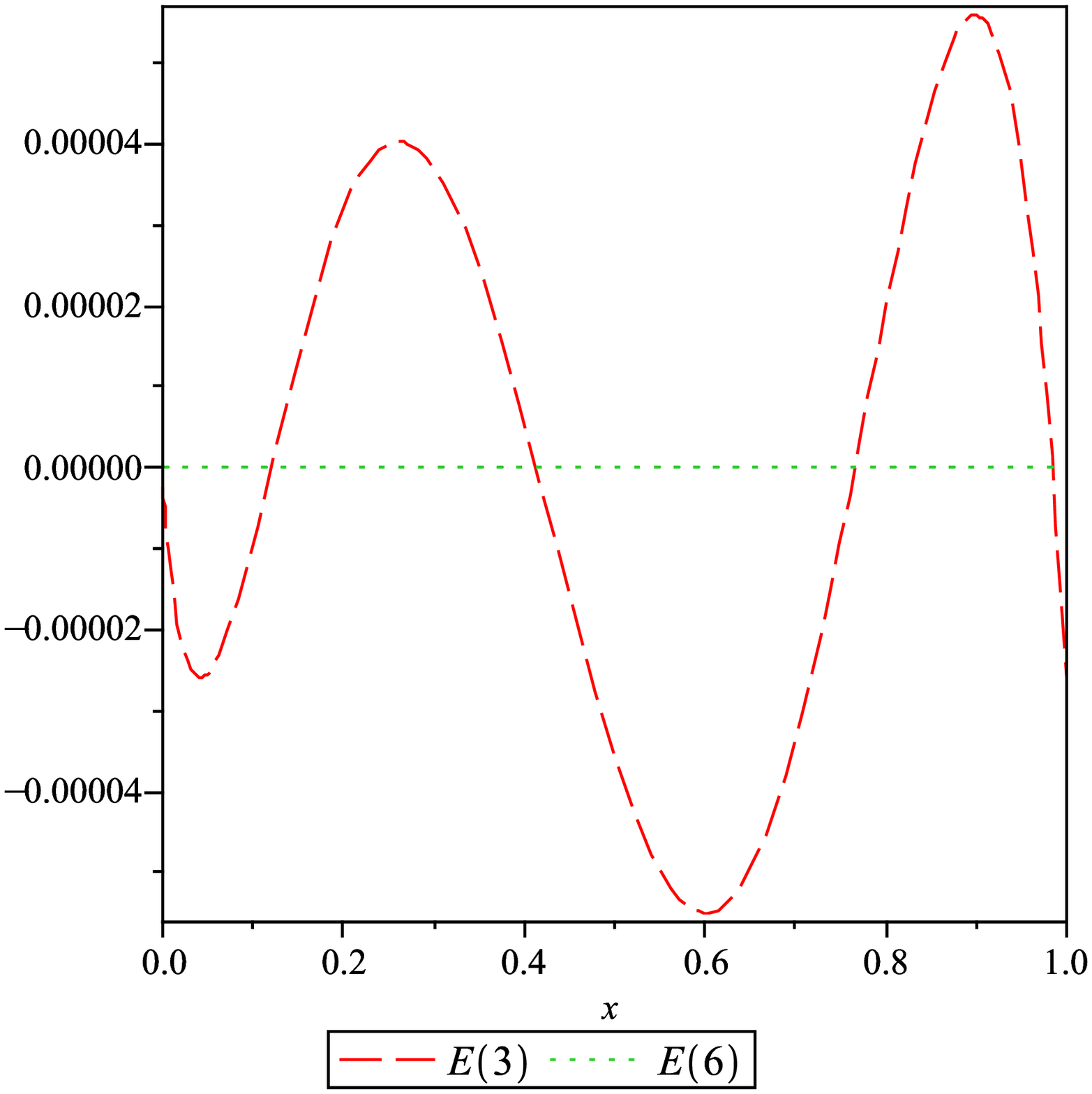}
  \caption{Comparison of exact solution and numerical solution and their errors for $n=3,\ 6$ and $\alpha=0.75$ and $ \epsilon=1$ in Example \ref{ex3}.}\label{Fig.3.}
\end{figure}

\begin{figure}
  \center\includegraphics[width=7cm]{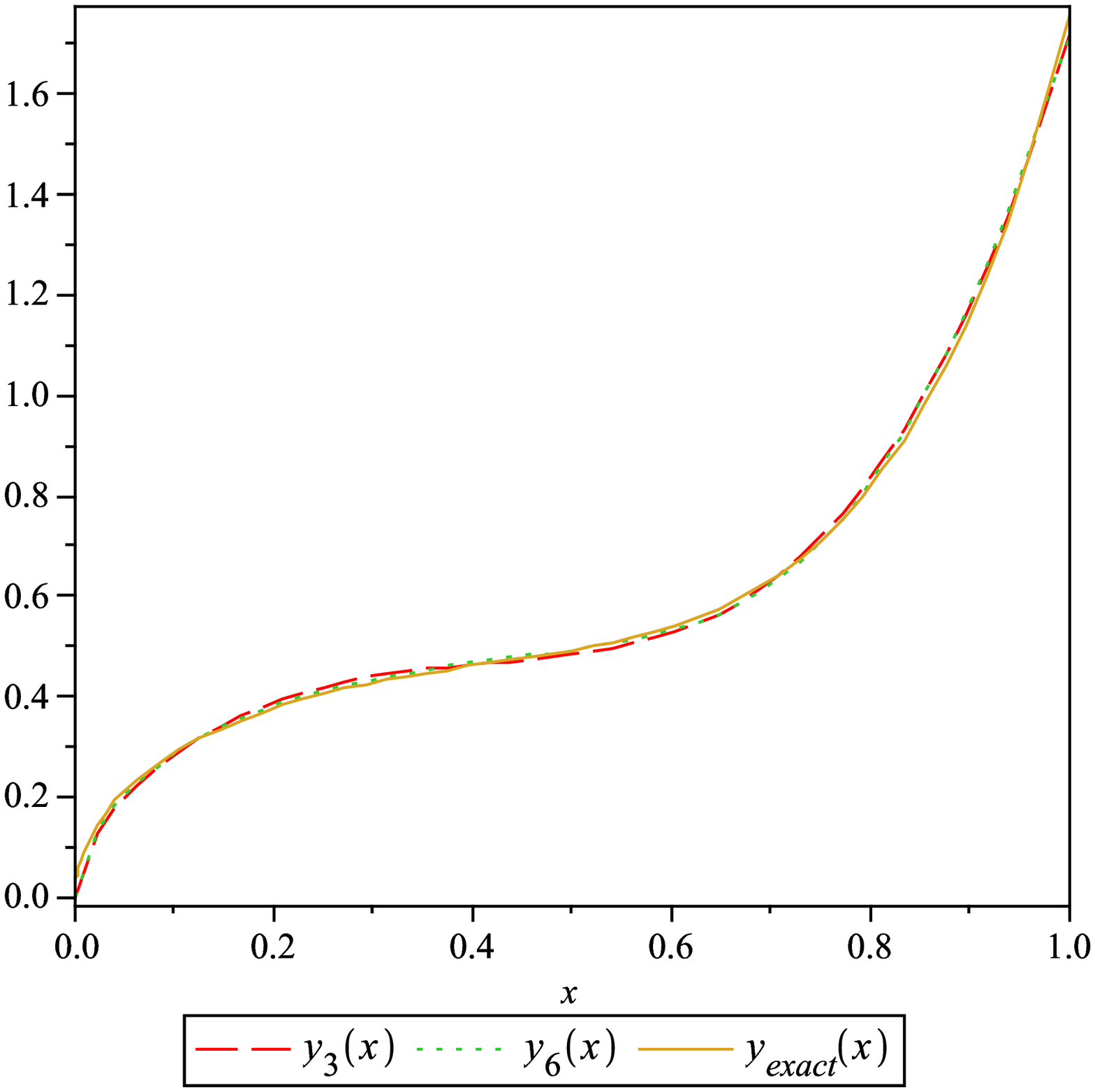}\includegraphics[width=7cm]{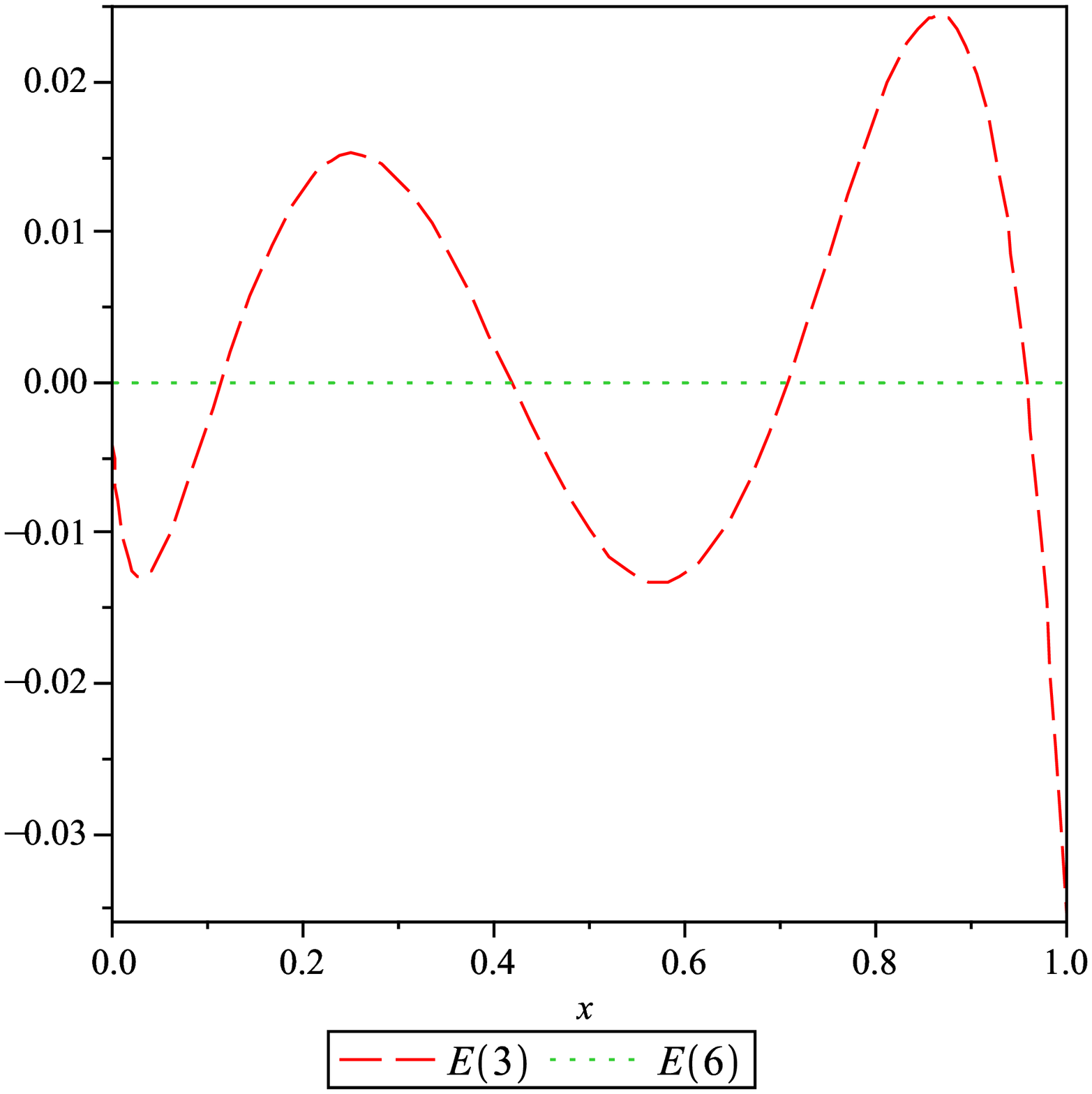}
  \caption{Comparison of exact solution and numerical solution and their errors for $n=3,\ 6$ and $\alpha=0.5,\ \beta=6$ and $ \epsilon=\nu=1$ in Example \ref{ex4}.}\label{Fig.4.}
\end{figure}

\begin{figure}
  \center\includegraphics[width=7cm]{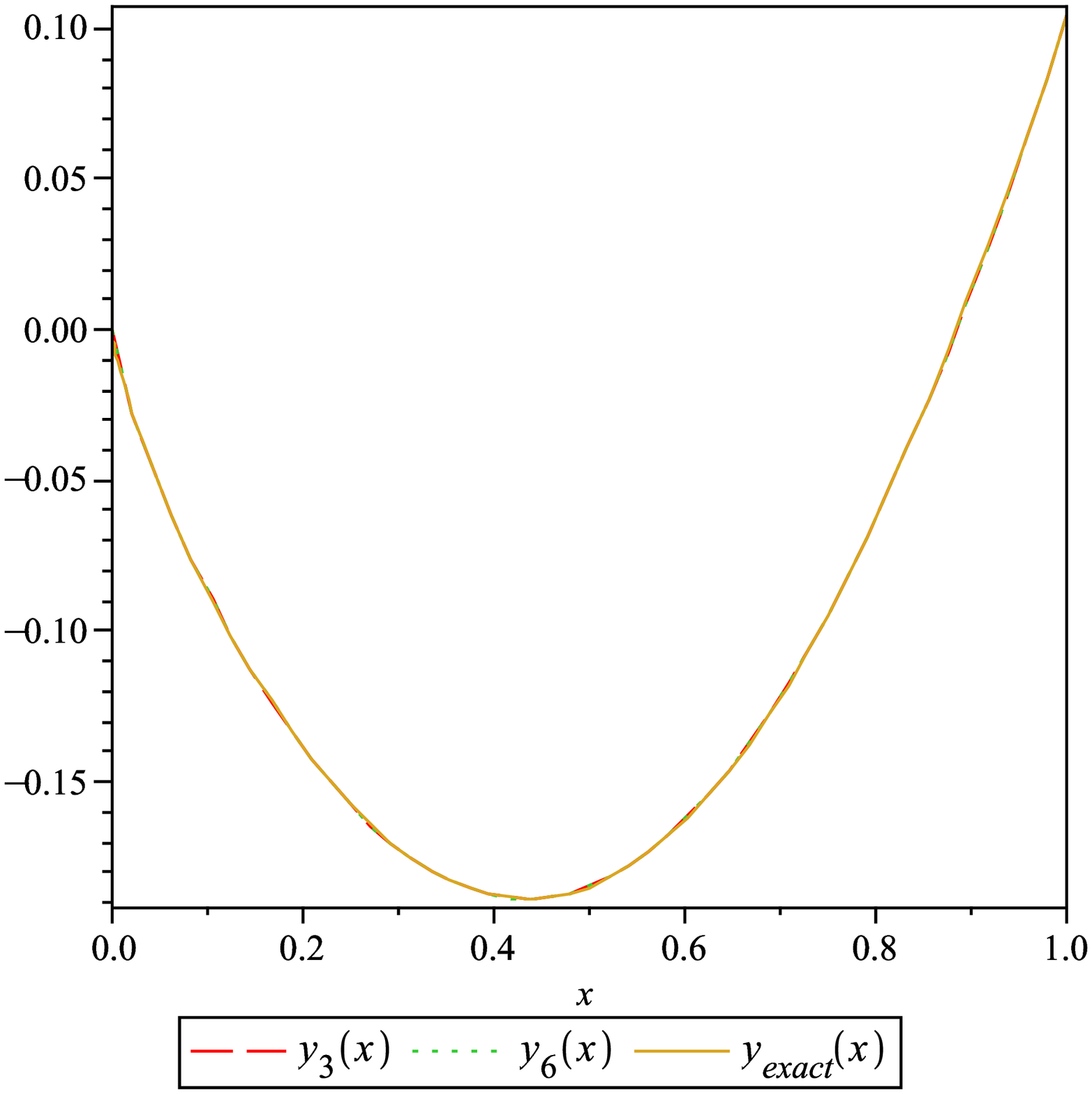}\includegraphics[width=7cm]{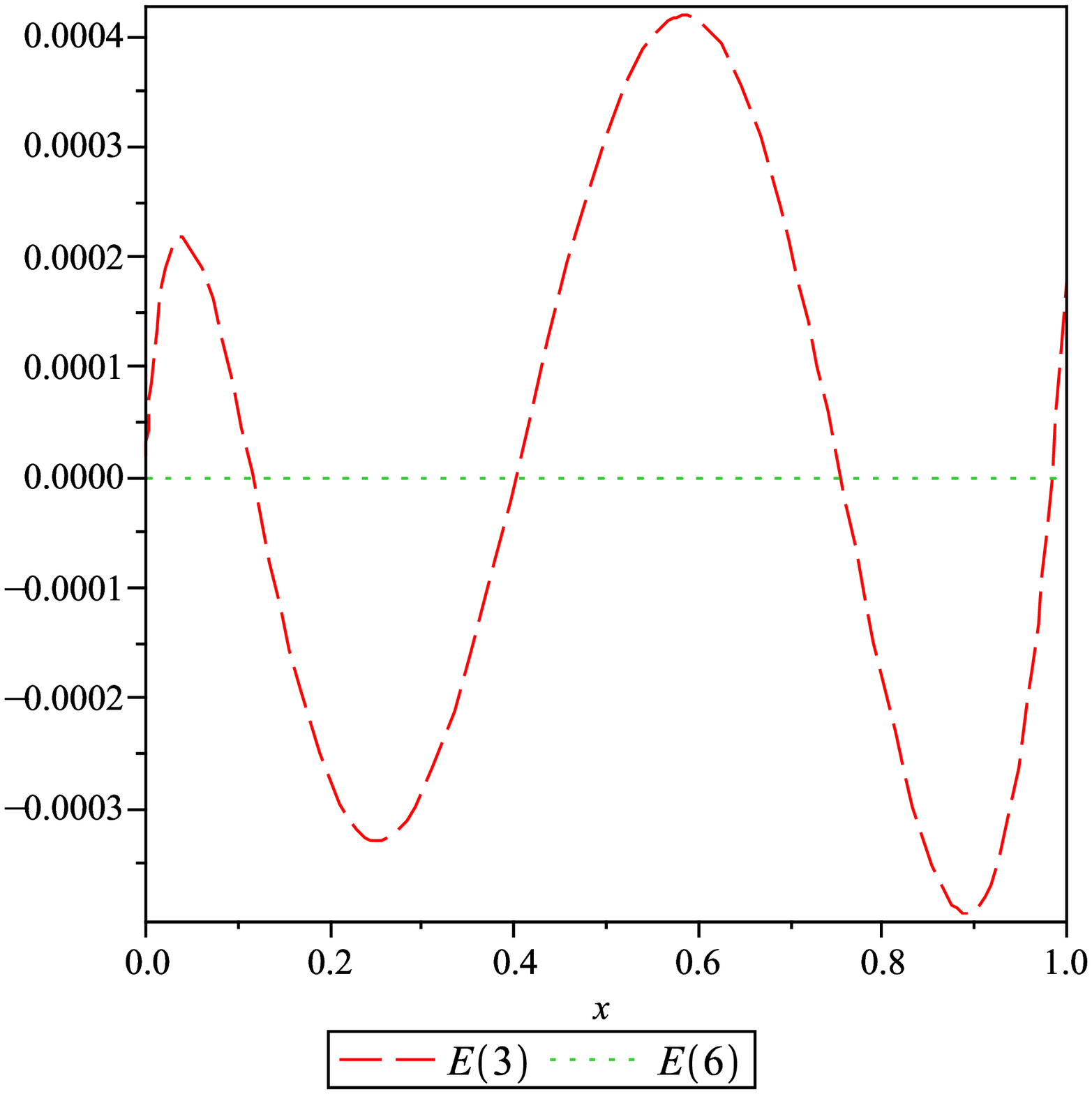}
  \caption{Comparison of exact solution and numerical solution  and their errors for $n=3,\ 6$ and $\alpha=0.75$ and $ \epsilon=0$ in Example \ref{ex5}.}\label{Fig.5.}
\end{figure}

\begin{Example}\label{ex1}
As first example, consider the fractional variational problem as in Theorem \ref{TeoExample} with $\displaystyle g(x)=h(x)=\frac{1}{1+x^\beta}$, then we have the following fractional variational problem:
\begin{equation}\label{probl1}
J[y]=\int_0^1 \left[\frac{1}{1+x^\beta}\ {}_{0}D_{x}^{\alpha}y(x)-({}_{0}I_{x}^{1-\alpha}y(x)+1)\frac{\beta x^{\beta-1}}{(1+x^\beta)^2}\right]^2dx \quad \rightarrow \quad \min,
\end{equation}
under the constraint
\begin{equation}\label{bcp1}
{}_{0}I_{x}^{1-\alpha}y(x)\Big|_{x=1}=\epsilon,
\end{equation}
In this case the exact solution is obtained from \eqref{exact} as:
\[
y_{exact}(x)=(\frac{1}{2}(1+\epsilon)-1)\left(\frac{\Gamma(\beta+2)}{\Gamma(\beta+\alpha+1)}x^{\beta+\alpha}+\frac{1}{\Gamma(\alpha+1)}x^{\alpha}\right)
+\frac{\Gamma(\beta+1)}{\Gamma(\alpha+\beta)}x^{\beta+\alpha-1}.
\]
\end{Example}
Comparison of exact solution and numerical solution for $n=3,\ 6$ and $\alpha=0.5$ for $\beta=5,\ \epsilon=1$ are shown in Fig.(\ref{Fig.1.}) (left). In Fig.(\ref{Fig.1.}) (right) error between exact solution and numerical solution $E(n)={y}_n(x)-y_{exact}(x)$ for $n=3,\ 6$ and $\alpha=0.5,\ \beta=5,\ \epsilon=1$ are shown.
\begin{Example}\label{ex2}
Consider now problem  of Theorem \ref{TeoExample} with $\displaystyle g(x)=h(x)=e^{-\nu x}$. Then, in this case,
\begin{equation}\label{probl2}
J[y]=\int_0^1 \left[e^{-\nu x} {}_{0}D_{x}^{\alpha}y(x)-\nu({}_{0}I_{x}^{1-\alpha}y(x)+1)e^{-\nu x}\right]^2dx \quad \rightarrow \quad \min,
\end{equation}
under the constrant
\begin{equation}\label{bcp2}
{}_{0}I_{x}^{1-\alpha}y(x)\Big|_{x=1}=\epsilon,
\end{equation}
The exact solution is
\[
y_{exact}(x)=(e^{-1}(1+\epsilon)-1)\nu^{-\alpha}\left(\sum_{k=0}^{\infty}\frac{(k+1)}{\Gamma(k+\alpha+1)}(\nu x)^{k+\alpha}\right)
+x^{\alpha-1}E_{1,\alpha}(\nu x)-\frac{x^{\alpha-1}}{\Gamma(\alpha)},
\]
where $E_{a,b}(x)$ is the Mittag-Leffler function of order $a$ and $b$ and defined as:
\[
E_{a,b}(x)=\sum_{k=0}^{\infty}\frac{x^{k}}{\Gamma(ak+b)}.
\]
\end{Example}
Exact solution and numerical solution for $n=3,\ 6$ and $\alpha=0.5$ and $\nu=1,\ \epsilon=-1$  are shown in Fig.(\ref{Fig.2.}) (left). Error between exact solution and numerical solution $E(n)={y}_n(x)-y_{exact}(x)$ for $n=3,\ 6$ and $\alpha=0.5$ and $\nu=1,\ \epsilon=-1$ are shown in Fig.(\ref{Fig.2.}) (right).
\begin{Example}\label{ex3}  For $\displaystyle g(x)=h(x)=\frac{1}{1+\sin(x)}$, the problem becomes
\begin{equation}\label{probl3}
J[y]=\int_0^1 \left[\frac{1}{1+\sin(x)}\  {}_{0}D_{x}^{\alpha}y(x)-({}_{0}I_{x}^{1-\alpha}y(x)+1)\frac{\cos(x)}{(1+\sin(x))^2}\right]^2dx \quad \rightarrow \quad \min,
\end{equation}
under the constraint
\begin{equation}\label{bcp3}
{}_{0}I_{x}^{1-\alpha}y(x)\Big|_{x=1}=\epsilon.
\end{equation}
The exact solution is
\[
y_{exact}(x)=A\left(\sum_{k=0}^{\infty}\frac{(2k+2)(-1)^k}{\Gamma(2k+\alpha+2)}x^{2k+\alpha+1}+\frac{x^{\alpha}}{\Gamma(1+\alpha)}\right)
+x^{\alpha}E_{2,\alpha+1}(-x^2),
\]
where $A=\displaystyle \frac{1}{1+\sin(1)}(\epsilon+1)-1$.
\end{Example}
Exact solution and numerical solution and their errors for $n=3,\ 6$ and $\alpha=0.75$ and $ \epsilon=1$  are shown in Fig.(\ref{Fig.3.}) (left) and (right), respectively.
\begin{Example}\label{ex4} Consider the fractional variational problem as in Theorem \ref{TeoExample} with $\displaystyle g(x)=\frac{1}{1+x^\beta}$ and $h(x)=e^{-\nu x}$, then we have the following fractional variational problem:
\begin{equation}\label{probl4}
J[y]=\int_0^1 \left[\frac{1}{1+x^\beta}\ {}_{0}D_{x}^{\alpha}y(x)-
\frac{\beta x^{\beta-1}}{(1+x^\beta)^2}\ {}_{0}I_{x}^{1-\alpha}y(x)-\nu e^{-\nu x} \right]^2dx \quad \rightarrow \quad \min,
\end{equation}
under the constraint
\begin{equation}\label{bcp4}
{}_{0}I_{x}^{1-\alpha}y(x)\Big|_{x=1}=\epsilon,
\end{equation}
In this case the exact solution is obtained from \eqref{exact} as:
\begin{eqnarray}
y_{exact}(x)&=&A\left(\frac{\Gamma(\beta+2)}{\Gamma(\beta+\alpha+1)}x^{\beta+\alpha}
+\frac{1}{\Gamma(\alpha+1)}x^{\alpha}\right)
+\left(\frac{1}{\Gamma(\alpha)}x^{\alpha-1}+\frac{\Gamma(\beta+1)}{\Gamma(\alpha+\beta)}x^{\beta+\alpha-1}\right)\nonumber\\
&-&\left(x^{\alpha-1}E_{1,\alpha}(-\nu x)+x^{\beta+\alpha-1} \sum_{k=0}^{\infty}\frac{\Gamma(k+\beta+1)(-1)^k}{\Gamma(k+\beta+\alpha)\Gamma(k+1)}(\nu x)^{k}\ \right),
\end{eqnarray}
where $\displaystyle A=-\frac{1}{2}+e^{-1}$.
\end{Example}
Comparison of exact solution and numerical solution and their errors for $n=3,\ 6$ and $\alpha=0.5,\ \beta=6$ and $ \epsilon=\nu=1$  are shown in Fig.(\ref{Fig.4.}) (left) and (right), respectively.
\begin{Example}\label{ex5} As the fifth  example,
consider the fractional variational problem as in Theorem \ref{TeoExample} with
$\displaystyle g(x)=\frac{1}{1+\sin(x)}$ and $h(x)=\cos(x)$, then we have the following fractional variational problem:
\begin{equation}\label{probl5}
J[y]=\int_0^1 \left[\frac{1}{1+\sin(x)}\  {}_{0}D_{x}^{\alpha}y(x)-\frac{\cos(x)}{(1+\sin(x))^2}\ {}_{0}I_{x}^{1-\alpha}y(x)-\sin(x)\right]^2dx \quad \rightarrow \quad \min,
\end{equation}
under the constraint
\begin{equation}\label{bcp5}
{}_{0}I_{x}^{1-\alpha}y(x)\Big|_{x=1}=\epsilon,
\end{equation}
In this case the exact solution is as:
\begin{eqnarray}
y_{exact}(x)&=&A\left(\sum_{k=0}^{\infty}\frac{(2k+2)(-1)^k}{\Gamma(2k+\alpha+2)}x^{2k+\alpha+1}+\frac{x^{\alpha}}{\Gamma(1+\alpha)}\right)
+\left(\frac{1}{\Gamma(\alpha)}x^{\alpha-1}+x^{\alpha}E_{2,1+\alpha}(-x^2)\right)\nonumber\\
&-&\left(x^{\alpha}E_{2,1+\alpha}(-4 x^2)+x^{\alpha-1}E_{2,\alpha}(-x^2) \right),
\end{eqnarray}
where $\displaystyle A=\cos(1)-1$.
\end{Example}
Comparison of exact solution and numerical solution  and their errors for $n=3,\ 6$ and $\alpha=0.75$ and $ \epsilon=0$  are shown in Fig.(\ref{Fig.5.}) (left) and in Fig.(\ref{Fig.5.}) (right), respectively.

\section{Conclusions}

In this paper we present a numerical treatment for fractional variational problems, by means of a decomposition formula based on Jacobi polynomials. Although we keep inside variational calculus, similar techniques could be used to solve fractional differential equations depending on fractional derivatives and fractional integrals of Riemann-Liouville type. In fact, it has already been done with success when in presence of Caputo fractional derivatives \cite{esmsha:11:pse}.

\section*{Acknowledgments}

The second author was supported by Portuguese funds through the CIDMA - Center for Research and Development in Mathematics and Applications, and the Portuguese Foundation for Science and Technology (“FCT–-Funda\c{c}\~{a}o para a Ci\^{e}ncia e a Tecnologia”), within project UID/MAT/04106/2013.

\end{document}